\documentclass[11pt]{amsart}
\usepackage{amssymb,latexsym,eucal,amsfonts,mathrsfs}
\newcommand{\iso}{\approx}

\newcommand{\ds}{\displaystyle}

\newcommand{\C}{{\mathbb{C}}}
\newcommand{\Z}{{\mathbb{Z}}}
\newcommand{\Q}{{\mathbb{Q}}}

\newcommand{\h}{{\mathfrak{h}}}
\newcommand{\g}{{\mathfrak{g}}}
\newcommand{\A}{\mathcal{I}}
\newcommand{\op}{\mathcal{A}}
\newcommand{\sh} {S(\h^*)}
\newcommand{\sind}[2]{S_{#1,#2}(w)}
\newcommand{\trans}[2]{t_{{#1}{#2}}}
\newcommand{\inv}[1]{#1^{-1}}
\newcommand{\sch}[1]{\mathfrak{S}_{#1}}
\newcommand{\ve}[3]{v_{[#1, #2, #3]}}

\begingroup
\newtheorem{bigthm}{Theorem}   
\newtheorem{thm}{Theorem}[section]   
\newtheorem{lemma}[thm]{Lemma}         
\newtheorem{prop}[thm]{Proposition}  
\newtheorem{defn}[thm]{Definition}   
   
\newtheorem{ex}[thm]{Example}        

\endgroup

\begin{document}


\title{A Pieri-type Formula for the Equivariant Cohomology of the Flag Manifold}
\author{Shawn Robinson}
\address{Department of Mathematics\\ 
  University of North Carolina at Chapel Hill\\
  Chapel Hill, NC}
\email{robinson@math.unc.edu}
\date{September 20, 2000}

\maketitle

\begin{abstract}
The classical Pieri formula is an explicit rule for determining the coefficients in the expansion 
$$\displaystyle s_{1^m} \cdot s_{\lambda} = \sum c_{1^m , \lambda}^{\mu} \ s_{\mu} \ ,$$
where $s_{\nu}$ is the Schur polynomial indexed by the partition $\nu$.  Since the Schur polynomials represent Schubert classes in the cohomology of the complex Grassmannian, this gives a partial description of the cup product in this cohomology.  Pieri's formula was generalized to the cohomology of the flag manifold in \cite{S}; and in the present work we generalize this formula to the $T$-equivariant cohomology of the flag variety. 

\end{abstract}

\section{Introduction} \label{S:intro}

Let $G$ be a semi-simple Lie group, $T$ be a maximal torus in $G$, $W$ be the corresponding Weyl group generated by simple reflections $\bigr\{ s_i \bigr\}_{1\leq i \leq l} ,$ $B$ be a Borel subgroup containing $T,$ and $\h$ be the Cartan subalgebra of $\g = Lie(G).$  Let $X = G/B.$  Denote by $ET \rightarrow BT$ the classifying bundle for the group $T.$  Since $T$ acts on $X$ we may form the fiber-product $ET \times_T X.$  The $T$-{\it equivariant cohomology} of $X$, denoted $H_T^*(X)$, is defined to be the singular cohomology $H^*(ET \times_T X;\Q).$  In \cite{KK}, B. Kostant and S. Kumar defined certain functions $\bigr\{ \xi^w :W \rightarrow S(\h^*) \bigr\}_{w \in W},$ and A. Arabia proved in \cite{A} that
$$ H_T^*(X) \iso \Lambda := \bigoplus_{w \in W}S(\h^*) \xi^w.$$
We call $\bigr\{ \xi^w \bigr\}_{w\in W}$ the {\it Schubert basis} of $\Lambda,$ and  the objects of study in this paper are the coefficients in the expansion
$$\xi^v \xi^w = \sum_{u \in W} p_{v, w}^u \xi^u.$$
We call these coefficients the {\it equivariant Schubert structure constants,} though these ``constants'' are actually polynomials in $\sh.$  In \cite{G}, W. Graham showed that these structure constants enjoy a certain positivity property, but it is an open problem to find a combinatorial expression for them that explicitly demonstrates this positivity.  The best previous result of this nature is a formula for multiplying any $\xi^w$ by a degree one Schubert class $\xi^{s_k}$ (\cite{KK}).

In this paper we consider only $G=SL_n(\C),$ in which case $W$ is the symmetric group $S_n$ and $X=G/B$ is the space of full flags of vector subspaces of $\C^n.$  Our main result is an {\it equivariant Pieri-type formula,} that is a formula for the structure constants in the expansion of the product $\xi^c \xi^w$, where $c=c[k,m]=s_{k-m+1} \dots s_{k-1}s_k.$   We recall a definition from \cite{M} (2.7.4) before stating our main result.

\begin{defn} \label{D:order}
Let $\trans{i}{j}$ denote the transposition $(i j)$, and let $\zeta = \trans{i_1}{q} \trans{i_2}{q} \dots \trans{i_p}{q} = (q\ i_p \dots i_2\  i_1),$  a $(p+1)$-cycle in $W$.  For any $w \in W$, we say $w\zeta$ is {\rm special $k$-superior to $w$ of degree $p$} if
\begin{enumerate}
\item $i_1 , \dots , i_p \leq k <q$,
\item $w(q) > w(i_1) > \dots >w(i_p)$, and
\item $l(w\zeta) = l(w) + p$.
\end{enumerate}

More generally, if $\zeta_1 , \zeta_2 , \dots , \zeta_r$ are pairwise disjoint cycles such that each $w\zeta_i$ is special $k$-superior to $w$ of degree $p_i$ and $\displaystyle \sum_{i=1}^r {p_i} = p$, then we say $w\zeta_1 \zeta_2 \dots \zeta_r$ is special $k$-superior to $w$ of degree $p$.  We denote the set of all such elements by $\sind{k}{p}$.
\end{defn}

\begin{bigthm} \label{T:pieri}
For any $w \in W$ and $c=c[k,m]$ as above, 
$$\xi^c \xi^w = \sum_{\substack{{0\leq p \leq m} \\ {u \in \sind{k}{p}}}} \xi^{c[k-p,m-p]}(\ve{u}{w}{k}) \ \xi^u,$$ 
where $\ve{u}{w}{k}$ is the element defined in \ref{D:assoc}. 
\end{bigthm}

In \cite{B}, S. Billey gave a compact expression for the values of the functions $\xi^w$ that clearly shows the positivity of the structure constants that appear in our Theorem \ref{T:pieri}.  In the special case $p=m$, our theorem implies the Pieri-type formula for Schubert polynomials first stated by A. Lascoux and M.-P. Schutzenberger in  \cite{LS1} and proved by F. Sottile in \cite{S}.  In \cite{LS2}, the authors introduced the {\it double Schubert polynomials} which represent the equivariant Schubert classes when $H^*_T(X)$ is identified with a quotient of a certain ring of polynomials.  When we take $n$ sufficiently large, our Pieri-type formula for $H^*_T(SL_n(\C)/B)$ implies a Pieri-type formula for double Schubert polynomials. 

In Section \ref{S:prelim} we cover preliminary material and fix our notation.  Section \ref{S:nilH} consists of a summary of relevant results from \cite{KK}. In Section \ref{S:result} we formulate our equivariant Pieri-type formula, state precisely the positivity property of the structure constants and verify it under the assumptions of Theorem \ref{T:pieri}, and compute an example.  We prove several lemmas concerning the set $\sind{k}{p}$ in section \ref{S:lemmas1}, and we establish some lemmas concerning the values of the $\xi^w$ in Section \ref{S:lemmas2}.  The proof of Theorem \ref{T:pieri} comprises Section \ref{S:proof}. 

We would like to thank especially S. Kumar for introducing us to this area of study and for many helpful discussions.  Thanks also to F. Sottile for referring us to \cite{M} which proved to be a very useful reference.

\section{Preliminaries} \label{S:prelim}
Throughout this paper (unless otherwise noted) $G$ denotes the group $SL_n(\C)$, and we choose the standard maximal torus, Borel subgroup, and set of simple roots.  Specifically, $T \subset G$ is the subgroup of diagonal matrices, $B$ is the subgroup of upper triangular matrices, $\h$ is the Cartan subalgebra corresponding to $T$, and $\{ L_1 ,\ L_2, \dots, L_n \} $ are the coordinates on $\h$ defined by 
$$L_i \Bigr( diag(t_1,\ t_2, \dots, t_n) \Bigr) = t_i.$$ 
The simple roots in $\h^*$ are $\bigr\{ \alpha_i = L_i - L_{i+1} \bigr\}_{1 \leq i \leq {n-1}},$ and the fundamental weights are $\ds \Bigr\{ \chi_i = \sum_{j=1}^{i} L_j \Bigr \}.$  The Weyl group $W = N(T)/T$ acts as the symmetric group permuting the subscripts of the coordinates $\{L_i\}$; more precisely, $W$ is generated by simple reflections $\{ s_i = s_{\alpha_i}\}_{1 \leq i \leq {n-1}}$ that act on $\h^*$ by
$$s_i \ L_j = 
\begin{cases}
  L_{i+1}, &\text{if $j=i$;}\\
  L_i , &\text{if $j=i+1$;}\\
  L_j, &\text{otherwise.}
\end{cases}$$
The group $W$ is characterized by the relations 
\begin{align}
&s_i^2 = 0, \notag \\
&s_i s_j = s_j s_i, \ \text{if $|i-j|> 1$, and}\notag \\
&s_i s_{i+1} s_i = s_{i+1} s_i s_{i+1}. \notag
\end{align}
  
The set of positive roots is $\Delta_+ =\{L_i - L_j \ | \ i <j \},$ and if $\beta = L_i - L_j$ is a positive root we denote the corresponding reflection in $W$ by $\trans{i}{j} = s_\beta.$  Whenever $\trans{i}{j}$ appears in this paper it is assumed that $i<j.$  We will often use $(i_1 \ i_2 \dots i_r)$ to denote the cycle $\zeta \in W$ such that $\zeta(i_j) = i_{j+1}$ for $1\leq j <r$, $\zeta (i_r) = i_1,$ and $\zeta(i)=i$ otherwise.  We say a decomposition $w= s_{i_1} s_{i_2} \dots s_{i_r} \in W$ is {\it reduced} if all other expressions of $w$ require at least $r$ simple reflections, in which case we say the {\it length} of $w$ is $r$ and write $l(w)=r.$  If $l(w\trans{i}{j})= l(w) +1$  we write $w \rightarrow w\trans{i}{j},$  and we write $w < u$ if there exists a chain
$w \rightarrow w\trans{i_1}{j_1} \rightarrow \dots \rightarrow w\trans{i_1}{j_1} \dots \trans{i_r}{j_r} = u.$

 For the flag manifold $X=G/B$ we have the famous Bruhat decomposition
$$X = \bigsqcup_{w \in W} BwB/B \text{ (disjoint union)},$$
where each $BwB/B$ is isomorphic to the affine space $\C^{l(w)}.$  Let $X_w$ be the closure of $BwB/B,$ the so-called {\it Schubert variety,} and let $[X_w]$ be the image of its fundamental homology class in $H_{2l(w)}(X; \Z).$  It is well known that these homology classes form a basis for $H_* (X; \Z).$  Let $\sch{w}$ be the class in $H^{2l(w)}(X; \Z)$ dual to $[X_w]$ with respect to the usual pairing of homology and cohomology -- we call these the {\it Schubert cohomology classes.}

\section{The $T$-equivariant cohomology of $G/B$} \label{S:nilH}

All the results in this section are from \cite{KK}, Sections 4 and 5, unless otherwise stated, and they hold for an arbitrary Ka\v{c}-Moody group $G,$ though we carry over the notation of the previous section.  The notational conventions here differ slightly from \cite{KK} and agree with \cite{K} and \cite{B}.  Let $S= S(\h^*),\ Q$ be its quotient field, and $\Omega$ be the $S$-algebra of all functions $W \rightarrow Q$ under point-wise addition and multiplication.  For each simple reflection $s_i$ we define a divided difference operator $\op_i$ on $\Omega$ by 
$$(\op_i f) (w) = \frac{f(ws_i) - f(w)}{w\alpha_i},$$
for any $f \in \Omega$ and any $w \in W.$  (These are analogous to the operators defined in \cite{BGG} and \cite{D}.)  Given a reduced decomposition $w=s_{i_1} \dots s_{i_l},$ we define $\op_w := \op_{i_1}\circ \dots \circ \op_{i_l},$ and this operator is independent of the choice of reduced decomposition of $w.$  
The setting for our Pieri-type formula is the ring
$$\Lambda := \bigr\{ f \in \Omega \ | \ f(w) \in S\  \text{and}\ (\op_w f)(e) \in S \ \text{for all}\  w\in W \bigr\}.$$

\begin{prop} 
The ring $\Lambda$ is a free $S$-algebra with basis $\{\xi^w \}_{w \in W},$ where the functions $\xi^w$ are characterized by:
\begin{enumerate}
\item $ \xi^w (v) = 0$ unless $w \leq v$;
\vspace{.05in}
\item $\ds \xi^w(w) = \prod_{\beta \in \Delta_+ \cap \inv{w}\Delta_{-}} \beta;$
\vspace{.05in}
\item $ \ds \op_i \xi^w =
  \begin{cases}
    \xi^{ws_i},  &\text{if $w > ws_i$}\\
    0, &\text{otherwise.}
    \end{cases}$
\end{enumerate}
\end{prop}
Furthermore, for any simple reflection $s_i$ the values of the function $\xi^{s_i}$ are
\begin{equation} \label{E:val}
\xi^{s_i} (v) = \chi_i - v\chi_i.
\end{equation}

The inclusion of the $T$-stable points of $X$ into $X$ induces a {\it localization map} on equivariant cohomologies:
$$ H_T^*(X) \rightarrow H^*_T(X^T).$$
If we equip $W$ with the discrete topology and let $T$ act trivially on $W,$ then $X^T \iso W$, and hence $H^*_T(X^T) \iso \Omega.$

\begin{thm} (\cite{A}, \cite{K})
The localization map induces an $S$-algebra isomorphism 
$$H_T^*(X) \stackrel{\sim}{\rightarrow} \Lambda.$$
Moreover, the map $\Lambda \rightarrow H^*(X)$, defined by $f\xi^w \mapsto f(0) \sch{w},$ for any $f \in S$, is an algebra homomorphism.
\end{thm}

In light of this theorem it is clear that formulas for the product in $\Lambda$ imply analogous formulas for the cup product in $H^*(X).$  We note here that in \cite{KK} a method for computing the equivariant structure constants in terms of the values $\xi^w(v)$ is described, and in \cite{B} the author gives a recursive method for computing the structure constants, also in terms of the values $\xi^w(v).$   In both cases though, the descriptions involve alternating sums and the cancellations that occur are not apparent.  In fact, our equivariant Pieri-type formula was found independently of either of these methods.      

\section{Equivariant Pieri-type formula} \label{S:result}

We begin this section by defining the element $\ve{u}{w}{k}$ that appears in Theorem \ref{T:pieri}. 
\begin{defn} \label{D:assoc}
Let $u \in \sind{k}{p}$ (see definition \ref{D:order}).   We define an associated element $\ve{u}{w}{k} \in W$ as follows.

{\rm Step 1:}  The disjoint cycles in the decomposition of $\inv{w}u$ are of the form $\zeta = \trans{i_1}{q} \dots \trans{i_r}{q} = (q \ i_r \dots i_1),$ with $i_1, \dots , i_r \leq k <q.$  Define the set of indices of $u > w$ to be 
$${\A} = \lbrace \  i \leq k \ |\ \inv{w}u (i) \neq i \rbrace.$$
In other words, if $u= w\trans{i_1}{q_1} \dots \trans{i_p}{q_p} \in \sind{k}{p}$ and $i_j \leq k <q_j$ for all $1 \leq j \leq p,$ then $\A = \{ i_1, \dots, i_p \}.$  The cardinality of $\A$ is $p$, and each member of $\A$ is at most $k.$

{\rm Step 2:}  Let $r_1$ be the greatest integer less than or equal to $k$ that is not in $\A ,$ and set $\lambda_1 =  \# \lbrace \ i \in {\A} \ | \ i < r_1 \rbrace.$ 
We define $\lambda_2, \dots \lambda_{k-p}$ recursively.  For each $1 \leq j \leq k-p-1,$ let $r_{j+1}$ be the greatest integer less than $r_j$ that is not in  $\A,$ and let $\lambda_{j+1} = \# \lbrace \ i \in {\A} \ | \ i <r_{j+1} \rbrace.$  The result is a partition, 
$$\lambda: \  p \geq \lambda_1 \geq \lambda_2 \geq \dots \geq \lambda_{k-p} \geq 0.$$     

{\rm Step 3:}  For each $j, \ 1 \leq j \leq k-p,$ let 
$$\pi_j = s_{\lambda_{k-p-j+1}+j-1} \dots s_{j+1}s_{j},$$
provided $\lambda_j \neq 0.$  If $\lambda_j = 0$, let $\pi_j = e.$  

{\rm Step 4:}  We define the associated element to be 
$$\ve{u}{w}{k} = w\pi_1 \pi_2 \dots \pi_{k-p}.$$
\end{defn}

\noindent With this definition in hand, we restate the main result of this paper.
\vspace{.1in}

\noindent {\bf Theorem A.}
{\it Let $w$ be an arbitrary element of $W$, and let $c=c[k,m].$  If $u$ is in $\sind{k}{p}$ for some $0 \leq p \leq m,$  then 
$$p_{c,w}^u = \xi^{c[k-p,m-p]}(\ve{u}{w}{k}).$$  
If $u$ does not belong to any such $\sind{k}{p}$, then $p_{c,w}^u = 0.$}  
\vspace{.1in}

Take as a basis for $S(\h^*)$ the simple roots $\alpha_1 , \dots , \alpha_{n-1},$ and consider the expansion
$$p_{v, w}^u = \sum_{I= (i_1, \dots, i_{n-1})} c_I \ \alpha_1^{i_1} \dots \alpha_{n-1}^{i_{n-1}}.$$  The following is proved in \cite{G} for an arbitrary Ka\v{c}-Moody group $G:$

\begin{prop}  
The coefficients $c_I$ in the expansion above are nonnegative integers for any $v, w, u \in W.$
\end{prop} 

To see that the coefficients in our Theorem \ref{T:pieri} have this positivity property we need a result from \cite{B}.  Given any $w\in W,$ fix a reduced decomposition $w=s_{i_1} \dots s_{i_l}.$  For each $1 \leq j \leq l,$ let $\beta_j = s_{i_1} \dots s_{i_{j-1}} \alpha_{i_j}.$  The following holds for any Ka\v{c}-Moody group $G$:
 
\begin{prop} \label{T:billey}
For any $v \in W$
$$\xi^v (w) = \sum \beta_{j_1} \dots \beta_{j_m},$$
where the sum is over all $1 \leq j_1 < \dots <j_m \leq l$ such that $s_{i_{j_1}} \dots s_{i_{j_m}}$ is a reduced decomposition of $v.$
\end{prop}

Together Theorem \ref{T:pieri} and Proposition \ref{T:billey} give an efficient method for computing the structure constants $p_{c, w}^u,$ as illustrated in the following example.

\begin{ex} {\rm
Let $G= SL_7(\C), \ c = c[4,2] = s_3 s_4, \ w= s_4 s_3 s_5 s_4,$ and $u=s_2 s_4 s_3 s_5 s_4 = w\trans{2}{5} \in \sind{4}{1}.$ Here the index set of $u>w$ is $\A =\{2\},$ so $\lambda$ is $(1,1,0),$ and the associated element is $\ve{u}{w}{4} = w \pi_1 \pi_2 \pi_3 = w (e)(s_2)(s_3) =  s_4 s_3 s_5 s_4 s_2 s_3.$  Thus
\begin{align}
p_{c,w}^u &=\xi^{s_3} ( s_4 s_3 s_5 s_4 s_2 s_3) \notag \\
&= s_4(\alpha_3) + s_4 s_3 s_5 s_4 s_2 (\alpha_3) \notag \\
&= \alpha_2 + 2\alpha_3 + 2\alpha_4 + \alpha_5. \notag
\end{align}} 
\end{ex}

\section{Lemmas regarding $\sind{k}{p}$} \label{S:lemmas1}

We begin this section by recalling the following lemma from \cite{BGG} (2.3{\it ii}).

\begin{lemma} \label{lemma0}
For any $w \in W$ and $i<j,$ $w < w\trans{i}{j}$ if and only if $w(j)>w(i).$

\end{lemma}

This implies that the chain of elements in definition \ref{D:order} is strictly increasing. More precisely, 
\begin{lemma} \label{lemma0.5}
If $w\trans{i_1}{q} \trans{i_2}{q} \dots \trans{i_p}{q} \in \sind{k}{p}$, then $w\trans{i_1}{q} \trans{i_2}{q} \dots \trans{i_r}{q} \in \sind{k}{r}$ for any $r \leq p$.
\end{lemma}

\begin{proof}
Since $w\trans{i_1}{q} \trans{i_2}{q} \dots \trans{i_p}{q} \in \sind{k}{p}$, we have $w(q)>w(i_1).$  By the preceeding lemma this is equivalent to $w < w\trans{i_1}{q}.$  Let $r$ be any index between $2$ and $p$.  It follows from $w\trans{i_1}{q} \dots \trans{i_{r-1}}{q} (q) = w(i_{r-1}) > w(i_{r}) = w\trans{i_1}{q} \dots \trans{i_{r-1}}{q} (i_r)$ that $w\trans{i_1}{q} \dots \trans{i_{r-1}}{q} < w\trans{i_1}{q} \dots \trans{i_{r}}{q},$  and by induction we have
$$w < w\trans{i_1}{q} < \dots < w\trans{i_1}{q} \dots \trans{i_r}{q} < \dots < w\trans{i_1}{q} \dots \trans{i_p}{q}.$$
Thus $l(w\trans{i_r}{q}) = l(w) + r.$ 
\end{proof}

The following four lemmas are needed for the proof of our main theorem.  They can be found in L. Manivel's proof of the Pieri-type formula for Schubert polynomials \cite{M}(2.7.5), but we provide more detailed proofs here.

\begin{lemma} \label{lemma1}
If $u \in \sind{k}{p} \setminus \sind{k-1}{p},$ then there is a unique index $q>k,$ such that $u\trans{k}{q} \in \sind{k-1}{p-1}.$  Furthermore, for any $q'<k,$ $u\trans{q'}{k} \notin \sind{k-1}{p-1}.$
\end{lemma}
\begin{proof}
The cycle decomposition of $w^{-1}u$ contains a cycle of the form $\zeta= \trans{i_1}{q}\dots \trans{i_a}{q} \dots \trans{i_r}{q}$, where $i_a =k$. To show that $u\trans{k}{q} \in \sind{k-1}{p-1}$ it suffices to show that $w\zeta \trans{k}{q} \in \sind{k-1}{r-1}.$ 
$$\zeta \trans{k}{q} = (\trans{i_1}{q} \dots \trans{i_{a-1}}{q})(\trans{i_{a+1}}{k} \dots \trans{i_r}{k}).$$
By Lemma \ref{lemma0} we have $w\zeta \trans{k}{q}< w\zeta,$ since $w\zeta (k)=w(i_{a-1}) > w(i_r)=w\zeta (q).$ By Lemma \ref{lemma0.5}, $w\trans{i_1}{q} \dots \trans{i_{a-1}}{q} \in \sind{k}{a-1},$ and by Lemma \ref{lemma0}  we have $l(w\trans{i_1}{q} \dots \trans{i_{a-1}}{q}\trans{i_{a+1}}{k} \dots \trans{i_r}{k} > l(w) + (a-1) + (r-a),$ since $w(i_a)>w(i_{a+1}) > \dots > w(i_r).$  Thus $l(w\zeta \trans{k}{q})=l(w)+r-1.$  It follows trivially from $w\zeta \in \sind{k}{r}$ that the remaining two conditions for $w\zeta \trans{k}{q} \in \sind{k-1}{r-1}$ are satisfied.

Now we show the uniqueness of this index $q.$  Pick any $q' > k, \ q' \neq q.$  The decomposition of $\inv{w}u$ may contain a cycle of the form $\zeta' = (q' \ j_s \dots j_1);$ then 
$$\zeta \zeta' \trans{k}{q'} = (q \ i_r \dots k \ j_s \dots j_1 \ q' \ i_{a-1} \dots i_1).$$Since this is an $(r+s+2)$-cycle, $w\zeta \zeta' \trans{k}{q'} \notin \sind{k-1}{r+s-1}.$  Hence $u\trans{k}{q'} \notin \sind{k-1}{p-1}.$  (Note that this argument holds when such a $\zeta'$ does not appear in $\inv{w}u.$  In this case delete all the $j$s.) 

Finally, we prove the second statement in the lemma.  Pick any $q'<k.$  In $\inv{w}u$ there may be a cycle of the form $\zeta' = \trans{j_1}{n} \dots \trans{j_b}{n} \dots \trans{j_s}{n},$ with $j_b = q'.$  For the second statement in the lemma it suffices to show that $w\zeta \zeta' \trans{q'}{k} \notin \sind{k-1}{r+s-1}.$ In cycle notation,
$$\zeta \zeta' \trans{q'}{k} = (q \ i_r \dots i_{a+1} \ k \ j_{b-1} \dots j_1 \ n \ j_s \dots j_{b+1} \ q' \ i_{a-1} \dots i_1).$$
This is an $(r+s+2)$-cycle, so $w \zeta \zeta' \trans{q'}{k} \notin \sind{k-1}{r+s-1}.$  Thus $u\trans{q'}{k} \notin \sind{k-1}{p-1}.$
\end{proof}

\begin{lemma} \label{lemma2}
If $u \in \sind{k}{p} \cap \sind{k-1}{p},$ then for any $q<k, \  u\trans{q}{k} \notin \sind{k-1}{p-1}.$  Also, for any $q' >k, \ u\trans{k}{q'} \notin \sind{k-1}{p-1}.$
\end{lemma}
\begin{proof}
First pick any $q<k.$  There are no cycles in $\inv{w}u$ that involve the index $k,$ but there may be a cycle of the form $\zeta = \trans{i_1}{n} \dots \trans{i_a}{n} \dots \trans{i_r}{n},$ where $i_a = q.$  In cycle notation 
$$\zeta \trans{q}{k} = (n \ i_r \dots i_{a+1} \ q \ k \ i_{a-1} \dots i_1).$$
Since this is an $(r+2)$-cycle, $w\zeta \trans{q}{k} \notin \sind{k-1}{r-1}.$  It follows that $u\trans{q}{k} \notin \sind{k-1}{p-1}.$

Now pick any $q'>k.$  In $\inv{w}u$ there may be a cycle of the form $\zeta = \trans{i_1}{q'} \dots \trans{i_r}{q'}.$  In this case $\zeta \trans{k}{q'} = (q' \ k\  i_r \dots i_1).$  Since this is an $(r+2)$-cycle, $w\zeta \trans{k}{q'} \notin \sind{k-1}{r-1},$ and it follows that $u\trans{k}{q'} \notin \sind{k-1}{p-1}.$
\end{proof}

\begin{lemma} \label{lemma3}
If $u \in \sind{k-1}{p} \setminus \sind{k}{p},$ then there is a unique $q<k,$ such that $u \trans{q}{k} \in \sind{k-1}{p-1}.$  Furthermore, for any $q'>k, \ u\trans{k}{q'} \notin \sind{k-1}{p-1}.$
\end{lemma}
\begin{proof}
The decomposition of $\inv{w}u$ contains a cycle of the form $\zeta = \trans{i_1}{k} \dots \trans{i_r}{k}.$  Let $q=i_r.$  By Lemma \ref{lemma0.5} we have $w\zeta \trans{q}{k} \in \sind{k-1}{r-1},$ and hence $u\trans{q}{k} \in \sind{k-1}{p-1}.$  The uniqueness of this index $q$ and the second statement in the lemma can be proved with arguments similar to those in the proof of Lemma \ref{lemma1}.
\end{proof}

\begin{lemma} \label{lemma4}
If $u \notin \sind{k}{p}$ and $u \notin \sind{k-1}{p},$ then there is a $q>k$ such that $u\trans{k}{q} \in \sind{k-1}{p-1}$ if and only if there is a $q' <k$ such that $u\trans{q'}{k} \in \sind{k-1}{p-1}.$ 
\end{lemma}
\begin{proof}
Suppose $u=u'\trans{k}{q}$ with $u' \in \sind{k-1}{p-1}$.  Since $u' \in \sind{k-1}{p-1}$ but $u'\trans{k}{q} \notin \sind{k}{p},$ the cycle decomposition of $\inv{w}u'$ must contain a cycle involving $k$ or $q.$  Assume both are present and write
$$\zeta_1 = \trans{i_1}{k} \dots \trans{i_r}{k} \ {\rm and }\ \zeta_2 = \trans{j_1}{q} \dots \trans{j_s}{q}.$$
Since $u' \rightarrow u'\trans{k}{q},$ we have $ u'(k) < u'(q),$  and it follows that $w(j_s) > w(i_r).$  Let $n$ be the least index, $1 \leq n \leq r,$ such that $w(j_s) >w(i_n).$
Then
\begin{align}
\zeta_1 \zeta_2 \trans{k}{q} &= (\trans{i_1}{k} \dots \trans{i_{n-1}}{k})(\trans{j_1}{q} \dots \trans{j_s}{q})\trans{i_n}{k} \trans{k}{q} (\trans{i_{n+1}}{q} \dots \trans{i_r}{q}) \notag \\
&=(\trans{i_1}{k} \dots \trans{i_{n-1}}{k})(\trans{j_1}{q} \dots \trans{j_s}{q})\trans{i_n}{q} \trans{i_n}{k}(\trans{i_{n+1}}{q} \dots \trans{i_r}{q}) \notag \\
&=(\trans{i_1}{k} \dots \trans{i_{n-1}}{k})(\trans{j_1}{q} \dots \trans{j_s}{q} \trans{i_n}{q} \dots \trans{i_r}{q})\trans{i_n}{k}. \notag
\end{align}

Let $\zeta_1' =\trans{i_1}{k} \dots \trans{i_{n-1}}{k}$ and $\zeta_2'=\trans{j_1}{q} \dots \trans{j_s}{q} \trans{i_n}{q} \dots \trans{i_r}{q}.$  Since $w\zeta_1'\zeta_2'(k) = w(i_{n-1}) > w(j_s) = w\zeta_1'\zeta_2'(i_n),$ we have $l(w\zeta_1'\zeta_2') \leq l(w) +r+s.$  Also, $l(w\zeta_1' \zeta_2') \geq l(w)+r+s,$ since $w(j_s) > w(i_n) > \dots > w(i_r).$  Hence $w\zeta_1'\zeta_2' \in \sind {k-1}{r+s},$ and it follows that $u\trans{i_n}{k} \in \sind{k-1}{p-1}.$  The converse is proved similarly. 
\end{proof}

\section{Lemmas regarding $\Lambda$} \label{S:lemmas2}

In this section we establish some lemmas concerning the values of certain $\xi^v$ on the associated element from Definition \ref{D:assoc}.  We begin with a lemma from \cite{KK}(4.24b) (see also \cite{K}(11.1.7b)) .

\begin{lemma} \label{cor1}
Let $w, \ v$ be arbitrary elements of $W$, and let $s_i$ be a simple reflection.  If $ws_i > w,$ then  
$$\xi^w (vs_i) = \xi^w (v).$$
More generally, if $l(wv) = l(w)+l(v),$ then 
$$\xi^w (uv) = \xi^w (u).$$
\end{lemma}

Next we prove four lemmas that will be used in conjunction with Lemmas \ref{lemma1} through \ref{lemma4}, respectively, in the proof of our main result. 
\begin{lemma} \label{lemmab1}
If $u\in \sind{k}{p}$ and $u\trans{k}{q} \in \sind{k-1}{p-1}$, then $\ve{u}{w}{k} = \ve{u\trans{k}{q}}{w}{k-1}.$
\end{lemma}
\begin{proof} 
Let $\A$ denote the set of indices for $u >w,$ and let $\A'$ denote the set of indices for $u\trans{k}{q} >w.$  Since $\A = \A' \cup \lbrace k\rbrace,$ all objects in the definition of $\ve{u}{w}{k}$ coincide with those in the definition of $\ve{u\trans{k}{q}}{w}{k-1}.$ 
\end{proof}

\begin{lemma} \label{lemmab2}
If $u \in \sind{k}{p} \cap \sind{k-1}{p}$, then 
\begin{enumerate}
\item $uL_k = \ve{u}{w}{k} L_{k-p}$, and 
\item $\xi^{c[k-p-1, n]} (\ve{u}{w}{k-1}) = \xi^{c[k-p-1, n]}(\ve{u}{w}{k}),$ for any $n \geq 1.$
\end{enumerate}
\end{lemma}

\begin{proof}
(1)   Let $\A$ denote the set of indices for $u >w.$  Since $u \in \sind{k}{p} \cap \sind{k-1}{p}, \ k \notin \A.$  It follows by definition that $r_1 = k, \ \lambda_1 = p,$ and $\pi_{k-p} = s_{k-1}s_{k-2} \dots s_{k-p}.$  We observe that the simple reflections $s_i$ occuring  in the reduced words for $\pi_1 , \dots , \pi_{k-p-1}$ all have $i < k-1.$  Hence
\begin{align}
\ve{u}{w}{k} L_{k-p} &= w\pi_{1} \dots \pi_{k-p} L_{k-p} \notag \\
&= w\pi_1 \dots \pi_{k-p-1} L_k \notag \\ 
&=w L_k \notag \\
&= u L_k. \notag
\end{align}

(2)  From the definition we see that $\ve{u}{w}{k-1} = \ve{u}{w}{k} \inv{\pi_{k-p}}.$  For $n \geq 1$, 
$ \ds c[k-p-1, n]= s_{k-p-n} \dots s_{k-p-1}$ and $\inv{\pi_{k-p}} = s_{k-p} \dots s_{k-1}$ have no simple reflections in common, so 
$$l(c[k-p-1, n] \inv{\pi_{k-p}}) = l(c[k-p-1, n]) + l(\inv{\pi_{k-p}}).$$  
By Lemma \ref{cor1}, 
\begin{align}
\xi^{c[k-p-1, n]} (\ve{u}{w}{k}) &= \xi^{c[k-p-1, n]} (\ve{u}{w}{k}\inv{\pi_{k-p}}) \notag \\
&= \xi^{c[k-p-1, n]} (\ve{u}{w}{k-1}). \notag  
\end{align}
\end{proof}

The proof of the next lemma is elementary, but it involves rather lengthy computations.  
  
\begin{lemma} \label{lemmab3}
If $u \in \sind{k-1}{p} \setminus \sind{k}{p}$ and $u\trans{q}{k} \in \sind{k-1}{p-1}$, then there is an element $\sigma \in W$ that has the following properties:
\begin{enumerate}
\item $u L_k = \sigma L_{k-p},$
\item $\xi^{c[k-p-1, N]}(\sigma) = \xi^{c[k-p-1, N]}(\ve{u}{w}{k-1})$ for all $N \geq 1,$ and 
\item $\xi^{c[k-p, m-p]}(\sigma) = \xi^{c[k-p, m-p]}(\ve{u\trans{q}{k}}{w}{k-1}).$
\end{enumerate}
\end{lemma}

\begin{proof}
Let $v = \ve{u}{w}{k-1},$ and let $\A, \ \lbrace r_i \rbrace , \ \lbrace \lambda_i \rbrace , \ \lbrace \pi_i \rbrace$ be the objects involved in the construction of this element.  Let $v' = \ve{u\trans{q}{k}}{w}{k-1},$ and let $\A', \ \lbrace r_i '\rbrace , \ \lbrace \lambda_i ' \rbrace , \ \lbrace \pi_i ' \rbrace$ be the objects involved in the construction of this element.  Since $\A =\A' \cup \lbrace q \rbrace,$ there is an $n$ such that $r_i = r_i '$ if $1 \leq i \leq n$, $r_{n+1}' = q,$ and $r_i = r_{i+1}'$ if $n+1 \leq i \leq k-p-1.$  If $\lambda_{n+1}' = 0,$ then $v=v',$ so we assume $\lambda_{n+1}' >0.$  Let
$$\sigma = v's_{k-p-n} \dots s_{k-p-2} s_{k-p-1}.$$

(1)  First observe that $u L_k = u\trans{q}{k} L_q = w L_q,$ and 
\begin{align}
\sigma L_{k-p} &= v's_{k-p-n} \dots s_{k-p-2} s_{k-p-1} L_{k-p} \notag \\
&=v' L_{k-p-n} \notag \\
&= w\pi_1' \dots \pi_{k-p}' L_{k-p-n} \notag \\
&= w \pi_1' \dots \pi_{k-p-n}' L_{k-p-n} \notag \\
&= w\pi_1' \dots \pi_{k-p-n-1}' L_{k-p-n+\lambda_{n+1}'} \notag \\
&= w L_{k-p-n+\lambda_{n+1}'}. \notag 
\end{align}
 
It remains to show that $q = k-p-n +\lambda_{n+1}',$ and this requires a small trick.  We recall that $\# \A = p$ and write
$$q = k-p- \lbrace (k-q-1)-(\#\A - \lambda_{n+1}' -1) \rbrace + \lambda_{n+1}'.$$
From $\A = \A' \cup \lbrace q \rbrace$ and $q=r_{n+1}'$ it follows that 
$$\#\A - \lambda_{n+1}' -1 = \#\lbrace i \in \A \ | \ i  >q \rbrace.$$
Also, $k-q-1 = \# \lbrace q+1, q+2, \dots , k-1 \rbrace,$ so 
\begin{align}
(k-q-1) - (\#\A - \lambda_{n+1}' -1) &= \# \lbrace i \notin \A \ | \ q< i \leq k-1 \rbrace \notag \\
&= n. \notag
\end{align}
Hence $q = k-p-n + \lambda_{n+1}'.$

(2)  Now we show by direct computation that $l(c[k-p-1, N] \inv{v}\sigma) = l(c[k-p-1, N])+ l(\inv{v}\sigma),$ for any $N \geq 1.$  It will then follow by lemma \ref{cor1} that 
$$\xi^{c[k-p-1, N]}(v)= \xi^{c[k-p-1, N]}(v \inv{v}\sigma) = \xi^{c[k-p-1, N]}(\sigma).$$
A rather laborious computation using the relationship between $\A$ and $\A'$ yields
$$\inv{v} \sigma = s_{k-p+\lambda_{n+1}'-1} \dots s_{k-p+1} s_{k-p},$$
and
$$c[k-p-1, N] \inv{v} \sigma = (s_{k-p-N} \dots s_{k-p-1})(s_{k-p+\lambda_{n+1}'-1} \dots s_{k-p}).$$
This is a reduced decomposition, so $l(c[k-p-1, N] \inv{v}\sigma) = l(c[k-p-1, N])+ l(\inv{v}\sigma).$

(3)  Finally, we observe that
$$c[k-p, N](s_{k-p-n} \dots s_{k-p-2} s_{k-p-1})  = (s_{k-p-N+1} \dots s_{k-p-1} s_{k-p})(s_{k-p-n} \dots s_{k-p-2} s_{k-p-1})$$
is a reduced word.  Hence 
$$l(c[k-p, N] s_{k-p-n} \dots s_{k-p-1}) = l(c[k-p, N]) +l( s_{k-p-n} \dots s_{k-p-1}),$$ 
and by Lemma \ref{cor1} 
$$\xi^{c[k-p, m-p]}(v') = \xi^{c[k-p, m-p]}(v's_{k-p-n} \dots s_{k-p-2} s_{k-p-1}).$$ 
\end{proof}

\begin{lemma} \label{lemmab4}
If $u \trans{k}{q} \in \sind{k-1}{p-1}$ and $u\trans{q'}{k} \in \sind{k-1}{p-1},$ but $u \notin \sind{k}{p}$, then $\ve{u\trans{k}{q}}{w}{k-1} = \ve{u\trans{q'}{k}}{w}{k-1}.$
\end{lemma}

\begin{proof}
The set {$\mathcal A$} of indices for $u\trans{k}{q} > w$ is also the set of indices for $u\trans{q'}{k} >w.$  $\square$  
\end{proof}

\section{Proof of the equivariant Pieri-type formula} \label{S:proof}

We begin this section by recalling a result from \cite{KK}(4.30) (see also \cite{K}(11.1.7{\it i})) that holds for any Ka\v{c}-Moody group $G.$ 
\begin{prop}  
$$\xi^{s_i} \xi^w = \xi^{s_i}(w) \ \xi^w + \sum_{w \rightarrow s_\beta w} \langle \chi_i , \inv{w}\beta \check{\ }\rangle \xi^{s_\beta v}$$
where $\beta \check{\ } \in \h$ is dual to $\beta \in \h^*$ with respect to the usual pairing $<\ ,\ >$ induced by the Killing form.  In the case $G= SL_n(C),$ this formula reads
\begin{equation} \label{chev}
\xi^{s_i}\xi^w = \xi^{s_i}(w)\  \xi^w + \sum_{\substack{w \rightarrow w\trans{i}{j} \\ i \leq k < j}} \xi^{w\trans{i}{j}}.
\end{equation}
\end{prop}
 
Let $e_m (X_k)$ denote the elementary symmetric function of degree $m$ in the variables $x_1,x_2, \dots , x_k$, and observe that 
\begin{align}
e_m(X_k) &= e_m(X_{k-1}) + x_k \ e_{m-1}(X_{k-1}). \notag \\
\intertext{In terms of the Schubert basis for the cohomology of $G/B$ this translates to} 
\sch{c[k,m]} &= \sch{c[k-1,m]} + (\sch{s_k}-\sch{s_{k-1}}) \ \sch{c[k-1,m-1]}. \notag 
\end{align}
The following key lemma is an analogous identity in the ring $\Lambda$ that follows directly from formula (\ref{chev}). 
\begin{lemma} \label{decomp1}
$$\xi^{c[k,m]} = \xi^{c[k-1,m]} + (\xi^{s_k} - \xi^{s_{k-1}} +L_{k-m+1} - L_k)\ \xi^{c[k-1,m-1]}$$
\end{lemma}
\begin{proof}
By Formula (\ref{chev}), 
\begin{align}
\xi^{s_k}  \xi^{c[k-1, m-1]} &= \xi^{s_k}(c[k-1, m-1])\ \xi^{c[k-1, m-1]}  +\xi^{c[k-2,m-2]s_{k}s_{k-1}} + \xi^{c[k,m]} \text{, and} \notag \\
\xi^{s_{k-1}} \xi^{c[k-1,m-1]} &= \xi^{s_{k-1}}(c[k-1, m-1])\ \xi^{c[k-1, m-1]} +\xi^{c[k-1,m]} + \xi^{c[k-2,m-2]s_k s_{k-1}} \text{, so} \notag\\
(\xi^{s_k}-\xi^{s_{k-1}}) \ \xi^{c[k-1,m-1]} &= (\xi^{s_k}(c[k-1, m-1]) -\xi^{s_{k-1}}(c[k-1, m-1]))\ \xi^{c[k-1, m-1]}+ \notag \\
&  \ + \xi^{c[k,m]} - \xi^{c[k-1, m]}. \notag 
\end{align}
By Formula (\ref{E:val}) from Section \ref{S:nilH}, 
\begin{align}
\xi^{s_k}(c[k-1, m-1]) -\xi^{s_{k-1}}(c[k-1, m-1]) &= \chi_k - c[k-1, m-1]\chi_k - \chi_{k-1} + c[k-1, m-1]\chi_{k-1}\notag \\
&= L_{k-m+1} - L_k. \notag 
\end{align}
\end{proof}
An obvious consequence of this lemma is 
\begin{equation} \label{decomp2}
\xi^{c[k,m]}(w) = \xi^{c[k-i,m]}(w) + (L_{k-m+1}-wL_k)\ \xi^{c[k-1,m-1]}(w).
\end{equation}
Also, this lemma and Formula (\ref{chev}) together imply
\begin{lemma} \label{coeff}
\begin{align}
p_{c[k,m], w}^u = & \ p_{c[k-1,m], w}^u + (L_{k-m+1} -uL_k) \ p_{c[k-1,m-1], w}^u + \notag \\ 
& + \sum_{u\trans{k}{q} \rightarrow u}p_{c[k-1,m-1], w}^{u\trans{k}{q}} - \sum_{u\trans{q}{k} \rightarrow u}p_{c[k-1,m-1], w}^{u\trans{q}{k}} \notag
\end{align}
\end{lemma}

\begin{proof}
Applying emma \ref{decomp1} we get
$$\xi^{c[k,m]} \xi^w = \xi^{c[k-1,m]}\xi^w + \Bigr\lbrace \xi^{c[k-1,m-1]}\xi^w \bigr(\xi^{s_k} -\xi^{s_{k-1}} +L_{k-m+1}-L_k\bigr) \Bigr\rbrace.$$
The term in braces may be rewritten as
$$\sum_v p_{c[k-1,m-1], w}^v \xi^v \bigr(\xi^{s_k} -\xi^{s_{k-1}} +L_{k-m+1}-L_k\bigr),$$
and we apply Formula (\ref{chev}) to obtain
$$\sum_v p_{c[k-1,m-1], w}^v \Bigr[\bigr(\xi^{s_k}(v)-\xi^{s_{k-1}}(v) +L_{k-m+1}-L_k\bigr) \xi^v + \sum_{\substack{v \rightarrow v\trans{i}{j} \\ i \leq k <j}}  \xi^{v\trans{i}{j}} \ - \sum_{\substack{v \rightarrow v\trans{i'}{j'} \\ i' \leq k-1 < j'}} \xi^{v\trans{i'}{j'}} \Bigr]=$$
$$= \sum_v p_{c[k-1,m-1], w}^v \Bigr[ \bigr(L_{k-m+1} -vL_k\bigr) \xi^v + \sum_{\substack{q>k \\ v \rightarrow v\trans{k}{q}}} \xi^{v\trans{k}{q}} - \sum_{\substack{q<k \\ v \rightarrow \trans{q}{k}}} \xi^{v\trans{q}{k}} \Bigr].$$
Thus
\begin{align}
\xi^{c[k,m]} \xi^w = & \ \xi^{c[k-1,m]}\xi^w + \sum_v  p_{c[k-1,m-1], w}^v(L_{k-m+1} -vL_k)\xi^v + \notag \\
 &+\sum_v \sum_{v \rightarrow v\trans{k}{q}} p_{c[k-1,m-1], w}^v \xi^{v\trans{k}{q}} -\sum_v \sum_{v \rightarrow v\trans{q}{k}} p_{c[k-1,m-1], w}^v \xi^{v\trans{q}{k}}, \notag
\end{align}
and gathering together the coefficients of $\xi^u$, we arrive at the conclusion of the lemma. 
\end{proof}

\noindent {\it Proof of Theorem \ref{T:pieri}}.  We induct on the index $k.$  If $k=1,$ then $m \leq 1$, and the theorem is exactly Formula (\ref{chev}).  We first consider the case $u \in \sind{k}{p},$ and we divide this into two subcases: $u \notin \sind{k-1}{p}$, and $u \in \sind{k-1}{p}$. 

  Suppose $u \in \sind{k}{p} \setminus \sind{k-1}{p}$.  Together Lemmas \ref{coeff} and \ref{lemma1} imply that 
$$p_{c,w}^u = p_{c[k-1, m], w}^u + p_{c[k-1, m-1], w}^{u} (L_{k-m+1}-uL_k) + p_{c[k-1, m-1], w}^{u\trans{k}{q}},$$
where $u\trans{k}{q}$ is the element of $\sind{k-1}{p-1}$ described in the proof of Lemma \ref{lemma1}.
By induction on $k$, the first two summands are zero, since $u \notin \sind{k-1}{p}$; and 
$$p_{c[k-1, m-1], w}^{u\trans{k}{q}} = \xi^{c[k-p,m-p]}(\ve{u\trans{k}{q}}{w}{k-1}).$$
By Lemma \ref{lemmab1}, $\ve{u\trans{k}{q}}{w}{k-1} = \ve{u}{w}{k}$, so the theorem is verified in this case.

Now suppose $u \in \sind{k}{p} \cap \sind{k-1}{p}$.  By Lemmas \ref{coeff} and \ref{lemma2}, 
$$p_{c,w}^u = p_{c[k-1, m], w}^u + p_{c[k-1, m-1], w}^{u} \bigr(L_{k-m+1}-uL_k\bigr),$$
and by induction,
$$p_{c,w}^u = \xi^{c[k-1-p, m-p]}(\ve{u}{w}{k-1}) + \bigr(L_{k-m+1}-uL_k\bigr) \xi^{c[k-1-p, m-1-p]}(\ve{u}{w}{k-1}).$$  
Applying Lemma \ref{lemmab2} and Formula (\ref{decomp2}), we get 
\begin{align}
p_{c,w}^u &= \xi^{c[k-1-p, m-p]}(\ve{u}{w}{k}) + (L_{k-m+1}-\ve{u}{w}{k}L_{k-p}) \xi^{c[k-1-p, m-1-p]}(\ve{u}{w}{k}) \notag \\
&= \xi^{c[k-p,m-p]}(\ve{u}{w}{k}). \notag
\end{align}

Next we consider the case $u \notin \sind{k}{p},$ and we examine separately two subcases:  $u \in \sind{k-1}{p}$, and $u \notin \sind{k-1}{p}$.  Suppose first $u \in \sind{k-1}{p}.$  By Lemmas \ref{lemma3} and \ref{coeff},
$$p_{c, w}^u = p_{c[k-1,m], w}^u + p_{c[k-1,m-1], w}^u \bigr(L_{k-m+1}-u L_k\bigr) - p_{c[k-1,m-1], w}^{u\trans{q}{k}},$$
where $u\trans{q}{k} \in \sind{k-1}{p-1}$ is as in the proof of Lemma \ref{lemma3}, and by induction,
\begin{align}
p_{c, w}^u = & \ \xi^{c[k-1-p,m-p]}(\ve{u}{w}{k-1}) + \bigr(L_{k-m+1}-u L_k\bigr)\xi^{c[k-1-p,m-1-p]}(\ve{u}{w}{k-1}) + \notag \\
&- \xi^{c[k-p, m-p]}(\ve{u\trans{q}{k}}{w}{k-1}). \notag
\end{align} 
Applying Lemma \ref{lemmab3} and Formula (\ref{decomp2}), we get
\begin{align}
p_{c, w}^u &= \xi^{c[k-1-p,m-p]}(\sigma) + \bigr(L_{k-m+1}-\sigma L_{k-p}\bigr)\xi^{c[k-1-p,m-1-p]}(\sigma) - \xi^{c[k-p, m-p]}(\sigma) \notag \\
&= \xi^{c[k-p,m-p]}(\sigma)  - \xi^{c[k-p, m-p]}(\sigma) = 0, \notag
\end{align}
so the theorem is verified in this case.

Finally, suppose $u \notin \sind{k-1}{p}$ and  there is a $q > k$ such that $u\trans{k}{q} \in \sind{k-1}{p-1}$.  By Lemma \ref{lemma4} there is a $q' < k$ such that $u\trans{q'}{k} \in \sind{k-1}{p-1},$ so it follows from Lemma \ref{coeff} that 
\begin{align}
p_{c, w}^u &= p_{c[k-1,m-1], w}^{u\trans{k}{q}} - p_{c[k-1,m-1], w}^{u\trans{q'}{k}} \notag \\
&= \xi^{c[k-p,m-p]}(\ve{u\trans{k}{q}}{w}{k-1}) - \xi^{c[k-p, m-p]}(\ve{u\trans{q'}{k}}{w}{k-1}). \notag
\end{align}
By Lemma \ref{lemmab4}, $\ve{u\trans{k}{q}}{w}{k-1}= \ve{u\trans{q'}{k}}{w}{k-1}$, so $p_{c, w}^u =0$ in this case, and the theorem is proved.

\end{document}